\documentclass[12pt]{article}
\newcommand{\reff}{\ref}\newcommand{\cit}[1]{\cite{#1}}\newcommand{\bibi}[1]{\bibitem{#1}}
\newcommand{\com}[1]{}
\usepackage{amssymb}
\usepackage{amsmath}
\usepackage{amsthm}
\usepackage[latin1]{inputenc}
\usepackage{graphicx}
\usepackage{color}
\usepackage{xypic}
\usepackage{epsf}

\newtheorem{lm}{{\bf Lemma}}
\newtheorem{theorem}[lm]{{\bf Theorem}}
\newtheorem{deff}[lm]{{\bf Definition}}

\newtheorem{cj}[lm]{{\bf Conjecture}}
\newtheorem{cor}[lm]{{\bf Corollary}}
\newtheorem{exemple}[lm]{{\bf Example}}

\newtheorem{prop}[lm]{{\bf Proposition}}

\theoremstyle{definition}
\newtheorem{rem}[lm]{{\bf Remark}}

\newcommand{\sep}{$\!${\rm.}\ }

\newcommand{\pr}{\noindent {\sl Proof\;\sep }}
\newcommand{\ep}{\hfill \framebox[2mm]{\ } \medskip}
\newcommand{\tq}{\ ;\ }
\newcommand{\be}{\begin{enumerate}}
\newcommand{\ee}{\end{enumerate}}
\newcommand{\noi}{\noindent}
\newcommand{\med}{\medskip}

\renewcommand{\a}{\alpha}

\renewcommand{\b}{\beta}
\newcommand{\cc}{{\bf C}}

\newcommand{\De}{\Delta}

\newcommand{\eps}{\varepsilon}

\newcommand{\K}{{\mathcal K}}

\newcommand{\q}{\theta}

\renewcommand{\S}{{\mathbb S}}

\newcommand{\N}{\mathbb{N}}
\newcommand{\Z}{\mathbb{Z}}

\newcommand{\R}{\mathbb{R}}

\newcommand{\card}{{\rm card}}

\newcommand{\ed}{\end{document}}

\title{Fair partitioning by straight lines}
\author{A.\ Fruchard and A.\ Magazinov\thanks{Supported by ERC Advanced Research Grant no. 267165 (DISCONV).}}
\date{\empty}
\begin{document}
\maketitle
\abstract{
A pizza is a pair of planar convex bodies $A\subseteq B$,
where $B$ represents the dough and $A$ the topping of the pizza. 
A partition of a pizza by straight lines is a succession of double operations:
a cut by a full straight line, followed by a Euclidean move of one of the
resulting pieces; then  the procedure is repeated.
The final partition is said to be fair if each resulting slice has 
the same amount of $A$ and the same amount of $B$.
This note proves that, given an  integer $n\geq2$, there exists a fair 
partition by straight lines of any pizza $(A,B)$ into $n$ parts if and only
if $n$ is even.
The proof uses the following result:
For any planar convex bodies $A, B$ with $A\subseteq B$, and any
$\a\in\,]0,\frac12[\,$, there exists an $\a$-section of $A$ which is a
$\b$-section of $B$ for some $\b\geq\a$. 
(An $\a$-section of $A$ is a straight line cutting $A$ into two parts, 
one of which has area $\a|A|$.)
The question remains open if the word ``planar'' is dropped.
}

\

{\noi\bf Keywords:} Convex body, alpha-section, fair partitioning.

\

{\noi\bf MSC Classification: 52A10, 52A38, 51M25, 51M04}

\

%
%
%
Let $\K$ denote the set of planar convex bodies, endowed with the
usual Hausdorf-Pompeiu metric. The area of $A\in\K$ is denoted by $|A|$
and its boundary is denoted by $\partial A$.
Following~\cit{cfv}, what we call a {\sl pizza} is a pair $(A,B)$ of two
nested planar convex bodies $A\subseteq B\subset\R^2$. 
We call $A$ the {\sl topping} and $B$ the {\sl dough}.
Given a pizza $(A,B)$ and an integer $n\geq2$,
a {\sl fair partition of $B$ in $n$ slices} is a family of $n$
internally disjoint convex subsets $B_1,\dots,B_n$  such that
$$
|B_1|=\dots=|B_n|~\mbox{ and }~|A\cap B_1|=|A\cap B_2|=\dots=|A\cap B_n|.
$$
For the sake of clarity, we call {\sl pieces}
the intermediate subsets and {\sl slices} the final ones.

There is a wide literature upon the  problem of fair partitioning a
convex body, see e.g.~\cit{kha}.
The expressions ``equipartition'' and ``balanced partition'' are also used.
If there is no other constraint than to obtain convex slices $B_i$,
then it has been proven recently that the answer is positive for 
all $n$, see e.g.~\cit{ka,s,so}.

In~\cit{bm} the authors use $k$-fans, which are half-lines 
starting from a common point. Since this process is very restrictive, 
the result is negative for $k\geq4$.

Other rules have also been considered. One can ask to have same perimeter 
and same area for each slice, see e.g.~\cit{bbs,bz}. In~\cit b, the author
uses only cuts by horizontal and vertical segments.
\med

In this note, we use a different cutting rule, which seems to be new:
Divide $B$ into two pieces with a straight cut. Each of the resulting
pieces is a convex body, their interiors are disjoint, and their union is $B$.
If $B$ is divided into $k$ pieces $B_{1},\dots,B_{k}$, choose one of
these pieces and divide it into two pieces with one straight cut.
After $n-1$ cuts, $B$ is divided into $n$ convex slices. We will refer to this rule
simply as to the {\sc cutting rule}, since no other rule is considered further in this note.

Our cutting rule is more restrictive than just partitioning $B$ into $n$ convex
bodies. For example, a non-degenerate $3$-fan partition cannot be obtained with
our rule.
However, for $k\geq4$, our rule becomes somewhat less restrictive
than a $k$-fan partition.
\med

%
%
%
Before going further, we need to introduce some notation.
The symbol $\S^1$ stands for the standard unit circle, $\S^1:=\R/(2\pi\Z)$,
endowed with its usual metric
$d(\q,\q')=\min\{|\tau-\tau'|\tq\tau\in\q,\;\tau'\in\q'\}$.
Given $\q\in\S^1$, let $\vec u(\q)$ denote the unit vector of direction $\q$,
$\vec u(\q)=(\cos\q,\sin\q)$ and let 
$\vec u\,'(\q)=\frac{d\vec u}{d\q}(\q)=(-\sin\q,\cos\q)$.

Given an oriented straight line $\De$ in the plane, $\De^+$ denotes the
{\sl closed} half-plane on the left bounded by $\De$, 
and $\De^-$ is the closed half-plane on the right.
We identify oriented straight lines with points of the cylinder 
$\cc=\S^1\times\R$, associating each pair $(\q,t)\in\cc$ to the line
oriented by $\vec u(\q)$ and passing at the signed distance $t$ from
the origin.
In other words, the half-plane $\De^+$ is given by 
$\De^+=\{x\in\R^2\tq\langle x,\vec u\,'(\q)\rangle\geq t\}$.
We endow $\cc$ with the natural distance 
$d\big((\q,t),(\q',t')\big)=\big(d(\q,\q')^2+|t-t'|^2)\big)^{1/2}$.
The reason to introduce the space $\cc$ is the following: several times throughout the paper we will say that some oriented
line moves continuously. The continuity will always refer to the topology of $\cc$. 

Given $\a\in\,]0,1[$ and $A\in\K$,
an {\sl $\a$-section of $A$} is an oriented line
$\De$ such that $|\De^-\cap A|=\a|A|$.
For all $\a\in\,]0,1[$ and all $\q\in[0,2\pi[$, there exists a unique
$\a$-section of $A$ of direction $\q$, denoted by $\De(\a,\q,A)$.
The line $\De(\a,\q, A)$, treated as a function, depends continuously on its three arguments.

Our first result has been conjectured in~\cit{cfv}.

\begin{theorem} \label{t1}
For any planar convex bodies $A, B$ with $A\subset B$, and any
$\a\in\,]0,\frac12[\,$, there exists an $\a$-section of $A$ which is a
$\b$-section of $B$ for some $\b\geq\a$.
\end{theorem}

\pr
By contradiction, if every $\a$-section $\De(\a,\q,A)$ of $A$ is a
$\b(\q)$-section of $B$ with $\b(\q)<\a$ then, by continuity of
$\q\mapsto\b(\q)$ and by compactness of $\S^1$, there exists $\eps>0$
such that, for all $\q\in\S^1$, $\b(\q)\leq\a-\eps$.

Choose an integer $n>\frac1\eps$.
Choose $x_0\in\partial A$ arbitrarily and,
for each positive integer $i\leq n$,
define $x_i$ recursively by $x_i\in\partial A$
and the oriented line $D_i=(x_{i-1}x_i)$ is an $\a$-section of $A$.

\begin{figure}[!h]
\centerline{\includegraphics[height=5cm]{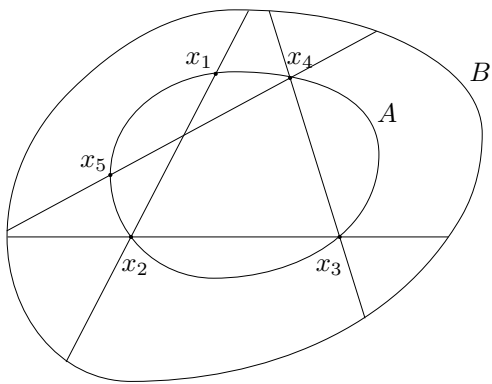}}
\vspace{3mm}
\caption{Construction of consecutive $\alpha$-sections}
\label{fig}
\end{figure}

We call a cap of $B$ the intersection $D_i^-\cap B$ and a cap of 
$A$ the intersection $D_i^-\cap A$. 
We thus have $n+1$ points on $\partial A$ and $n$ caps of $A$, resp. $B$.
For each $x\in B$, let $K(x)$ be the number of caps of $B$ which cover $x$,
i.e. $K(x)=\card\{i\in\{1,\dots,n\}\tq x\in D_i^-\}$.
Let $k=\min\{K(x)\tq x\in\partial A\}$; it is the number of complete tours
made by $x_0,\dots,x_n$.
Observe that, for all $x\in\partial A$ we have $k\leq K(x)\leq k+1$
and that, for each $0<m\leq n$,
the arc $\widehat{x_{m-1}x_m}=\partial A\cap D_m^-$ contains
$K(x_m)+1$ or $K(x_m)+2$ points among $x_0,\dots,x_n$
(including $x_{m-1}$ and $x_m$).
This comes from the fact that, if $x_i$ is on the arc $\widehat{x_{m-1}x_m}$,
then $x_m$ is on the arc $\widehat{x_{i}x_{i+1}}$.

We claim that $K(x)\leq k+1$ for all $x\in A$.
Indeed, let $x\in A$ and consider a cap $A_m$ containing $x$
(if $x$ belongs to no cap, there is nothing to prove).
If another cap $A_i$ contains $x$, then $x_{i-1}$ or $x_i$ must belong 
to the open arc $\widehat{x_{m-1}x_m}\setminus\{x_{m-1},x_m\}$.
Now for each of the points $x_j$ on this open arc,
at most one of the caps $A_j$ or $A_{j+1}$
can contain $x$, hence $x$ is in at most $k+1$ caps.
It follows that the sum of areas of all caps
$D_i^-\cap A$ is at most $(k+1)|A|$, hence $n\a\leq k+1$.

Also, for all $x\in B\setminus A$, we have $K(x)\geq k$. Hence 
$\sum_{i=1}^n|D_i^-\cap(B\setminus A)|\geq k|B\setminus A|$. Thus there exists 
$i_0$ such that $|D_{i_0}^-\cap(B\setminus A)|
\geq\frac kn|B\setminus A|\geq(\a-\frac1n)|B\setminus A|$.
It follows that 
$|D_{i_0}^-\cap B|=|D_{i_0}^-\cap A|+|D_{i_0}^-\cap(B\setminus A)|
\geq\a|A|+(\a-\frac1n)|B\setminus A|\geq(\a-\frac1n)|B|$,
i.e. $D_{i_0}$ is a $\b$-section of $B$,
with $\b\geq(\a-\frac1n)>\a-\eps$,
a contradiction.
\ep

\begin{rem}
A question whether Theorem~\ref{t1} extends to an arbitrary dimension remains open. However, one can show that for every $d > 2$
there exists a constant $\alpha_0(d) > 0$ such that the $d$-dimensional analogue of Theorem~\ref{t1} holds for all $\alpha \in ]0, \alpha_0[$.
The idea is similar to the 2-dimensional proof, but instead of an $n$-fold covering of $\partial A$ by caps we use a 1-fold covering,
namely, the {\it economic cap covering}, defined, for example, in~\cite{bar_rp}. However, this method is not very efficient, 
giving only a very small value of $\alpha_0$. Hence we leave the details of the proof to the reader.  
\end{rem}

Another equivalent formulation of Theorem~\reff{t1}, which will be more
convenient, is as follows. The proof of the equivalence is easy and left
to the reader.

\begin{cor}\label{c1}
For any planar convex bodies $A, B$ with $A\subset B$, and any
$\a\in\,]0,\frac12[\,$, there exists an $\a$-section of $B$ which is a
$\b$-section of $A$ for some $\b\leq\a$.
\end{cor}

Our next result, Theorem~\ref{t2}, concerns a fair pizza partition problem using the {\sc cutting rule}.
It has been already mentioned in~\cit{cfv} as a consequence of Theorem~\ref{t1}, but without a proof of implication.
Here we give a proof, and thus confirm the result. 
%
\begin{theorem}\label{t2}
Let $n$ be a positive integer. Then
\begin{enumerate}
	\item[\rm 1.] If $n$ is even, then for every pair $A\subseteq B$ of nested planar convex bodies there exists a fair partition obeying the
				{\sc cutting rule}.
	\item[\rm 2.] If $n$ is odd, then for some pairs $A\subseteq B$ such a partition may not exist.
\end{enumerate}
\end{theorem}

\pr
It is easy to check that two concentric disks $A$ and $B$ cannot be divided 
in a fair way into an odd number of slices: 
The first cut divides the pizza in two pieces,
containing $k$, resp. $l$ final slices,
with $k+l=n$ odd, hence $k\neq l$, and the smaller piece will not have enough
topping. 

To construct a fair partition for all even $n$, we proceed by induction.
 
For $n=2$, this follows from the intermediate value theorem.
Given $n\in\N$, $n$ even $\geq4$, and a pair of nested convex bodies
$A\subseteq B$, assume that a fair partitioning exists for any pair of nested
convex bodies and any even integer $i < n$.

If $n = 4k$, then the intermediate value theorem yields a fair
partitioning of two equal halves, and, by induction hypothesis, each of these halves
admits a fair partitioning in $2k$ slices. 

Let $n = 4k + 2$. Then we set $\a=\frac{2k}{4k+2}$ and consider two subcases.

\noindent {\bf 1.} Suppose that we can cut $B$ into two convex pieces $B_1$ and $B_2$
of areas $|B_1|=\a|B|$, $|B_2| = (1 - \a) |B|$ so that $|A\cap B_1|=\a|A|$, $|A\cap B_2|=(1-\a)|A|$.

Then, by induction, $B_1$ and $B_2$ have both
a fair partitioning in $2k$, resp. $2k+2$, slices, and this gives a fair
partitioning of $B$ in $n$ slices.

\noindent {\bf 2.} If we are not in subcase {\bf 1} then no $\a$-section of $B$ contains
an $\a$-portion of $A$. Then from Corollary~\reff{c1} it follows
that each $\alpha$-section of $B$ (with this $\a$) is a $\beta$-section of
$A$ for some $\beta < \alpha$.

Cut $B$ into two fair halves $B'$ and $B''$. We claim that there is a cut of $B'$
(and, similarly, of $B''$) such that it produces a slice of area $\frac{1}{n} |B|$ with
the topping part of area $\frac{1}{n} |A|$ (i.e., a fair slice).

Consider a piece $C_1 \subset B'$ between two parallel lines, one of which is the initial cut,
and the other one is chosen so that $|C_1| = \frac{1}{n} |B|$.
By the construction, $B' \setminus C_1$ is an $\a$-section of $B$, so $|A \cap (B' \setminus C_1)| < \a |A|$
and hence $|A \cap C_1| > \left( \frac 12 - \a \right) |A| = \frac 1n |A|$. 

On the other hand, by Corollary~\reff{c1},
there exists a $\frac 2n$-section of $B'$, which is at most $\frac 2n$-section of $A \cap B'$. If $C_2$ is the piece of $B'$
obtained by that section, then $|C_2| = \frac 1n |B|$, and $|A \cap C_2| \leq \frac 1n |A|$. 

Using the intermediate value theorem for $\frac 2n$-sections of $B'$,
we obtain that there is a slice $C$, which is cut from $B'$ by a single line, such that
$|C| = \frac 1n |B|$, and $|A \cap C| = \frac 1n |A|$.

By induction hypothesis, the piece $B' \setminus C$ admits a fair partition into $2k$ slices. As a result, there is a fair partition of $B'$
into $2k + 1$ slices. The same can be done for $B''$, yielding a fair partition of the whole pizza.

\ep

\

{\noi\bf Acknowledgements.}
The authors thank Maud Chavent from Plougonver who asked the question of partition,
and Nicolas Chevallier, Costin V\^\i{}lcu,
Imre B\'ar\'any, and Attila P\'or for fruitful discussions.

\vspace{1cm}

{\small 
\noi Addresses of the authors:
\med

\noi Augustin Fruchard\\
Laboratoire de Math\`ematiques, Informatique et Applications\\
Facult\`e des Sciences et Techniques\\
Universit\`e de Haute Alsace\\
2 rue des Fr\'eres Lumi\'ere\\
68093 Mulhouse cedex, FRANCE\med\\
E-mail: {\tt Augustin.Fruchard@uha.fr}

\

\noi Alexander Magazinov\\
Steklov Mathematical Institute\\ 
8 Gubkina Str.\\
Moscow 119991, Russia\\
E-mail: {\tt magazinov-al@yandex.ru}

\end{document}